\documentclass[11pt
]{elsarticle}
\usepackage{float}
\usepackage{amssymb}
\usepackage{amsmath}
\usepackage{amscd}
\usepackage{amsfonts}
\usepackage{amsthm}
\usepackage{graphicx}
\usepackage{mathtools}
\usepackage{color}
\usepackage[colorlinks, linkcolor=red, citecolor=green,filecolor=magenta,urlcolor=cyan]{hyperref}
\newfont{\Blackboard}{msbm10 scaled 1200}

\newfont{\roma}{cmr10 scaled 1200}

\newtheorem{Theorem}{{}\hskip\parindent Theorem}[section]
\newtheorem{lemma}{{}\hskip\parindent Lemma}[section]

\newtheorem{Example}{{}\hskip\parindent Example}[section]

\newtheorem{Definition}{{}\hskip\parindent Definition}[section]
\newtheorem{Remark}{{}\hskip\parindent Remark}[section]
\newtheorem{Proposition}{{}\hskip\parindent Proposition}[section]

\def\be{\begin{equation}}
	\def\ee{\end{equation}}
\def\beq{\arraycolsep=1.5pt\begin{eqnarray}}
	\def\eeq{\end{eqnarray}}

\begin{document}
	
	\begin{frontmatter}
	
	\title{ Lyapunov Stability for nonautonomous systems on Manifolds
	}

\date{}
\author[SWJTU]{Li Deng\corref{cor}}
\ead{dengli@swjtu.edu.cn}
\address[SWJTU]{School of Mathematics, Southwest Jiaotong University, Chengdu 611756, Sichuan Province, China}

\author[SWJTU]{Xin Li}
\ead{2195337638@qq.com}
%\address[Deng]{School of Mathematics,   Southwest Jiaotong University, Chengdu 611756, China. E-mail}
%address: 2195337638@qq.com}

\cortext[cor]{Corresponding author.
\\  This work is  supported by the National Science
Foundation of China under grant 12371451, and the Natural Science Foundation of Sichuan province under grants 2025ZNSFSC0077. 
}

	\begin{abstract}
This paper studies the uniformly asymptotic stability of nonautonomous systems on Riemannian manifolds. We establish corresponding Lyapunov-type theorems (Theorems~\ref{the1} and~\ref{the2}), extending classical Euclidean results (e.g., \cite[Theorems~4.9 and 4.10]{jind}) to curved spaces.
Our main contributions are: (i) an estimate for the domain of attraction linked to the equilibrium point's injectivity radius, where, under suitable conditions, this radius can be bounded using the sectional curvature (Proposition~\ref{cur1});
(ii) a demonstration that this estimate depends on the choice of the Riemannian metric (Examples~\ref{exa1} and ~\ref{exa2} and Remark \ref{r2d4}); 
 and (iii) a refined estimate compared to the Euclidean case, as detailed in  item (6) of Remark \ref{rt21}  and in Remark  \ref{r233}.
	\end{abstract}
	
	\begin{keyword}
	nonautonomous systems, 	uniformly asymptotic stability , Lyapunov-type theorem, Riemannian manifold.
		
	\end{keyword}
	
	\end{frontmatter}

%	\maketitle
\setcounter{equation}{0}		
	\section{Introduction}
\def\theequation{1.\arabic{equation}}
Let $n\in \mathbb{N}^+$, and let $M$ be a complete,  $n$-dimensional Riemannian manifold with the metric $g$. For $x\in M$, denote by  $T_{x}M$  the tangent space of $M$ at $x$,  by $T_{x}^{*}M$  the cotangent space, by $\langle\cdot,\cdot\rangle$ the inner product on $T_xM$, by $|\cdot|$ the norm induced by $\langle\cdot,\cdot\rangle$, and by $TM=\bigcup_{x\in M}T_{x}M$  the tangent bundle of $M$.  We denote by $\rho(\cdot,\cdot) $ the distance function on $ M $ with respect to $g$.    For standard definitions, see \cite[Sections 0.2, 1.2, 7.2]{c}.
	 	
	 Let $D \subseteq M$ be an open set. Let $f:[0,+\infty) \times D \to TM$ be a function (satisfying suitable assumptions to be given later). Let us consider the following ordinary differential equation
	 \begin{equation}\label{zuic}
	 	\dot{x}(t)=f(t,x(t)),\quad\forall\, t\geq t_0\geq  0,		
	 		 \end{equation}
	 where $\dot{x}(t)=\frac{d}{dt}x(t)$ for $t\geq 0$.   A point $x^{\ast}\in D$ is  referred as an equilibrium point of (\ref{zuic}), if $f(t,x^{\ast})=0,\;\forall\, t\geq 0$. 
	 
	 Regarding system (\ref{zuic}), we first introduce some fundamental concepts concerning Lyapunov stability.

\begin{Definition}\label{dsga}
	   The equilibrium point $x^*$ of system (\ref{zuic}) is uniformly asymptotically stable if there exists a neighborhood $\mathcal O$ of $x^*$, independent of $t_0\in[0,+\infty)$, such that for all $x(t_0)\in \mathcal O$, the solution  $x(\cdot)$ of the system (\ref{zuic}) satisfies: $ x(t)$ tends to $x^*$ as $t\to+\infty$, uniformly in $t_0$. Moreover, if there exists a subset  $\widetilde{D} \subseteq D$ such that for all $x(t_0)\in \widetilde D$, the corresponding solution $x(\cdot)$ satisfies  $x(t)$ tends to $x^*$  as $t\to+\infty$,  uniformly in $t_0$, then $x^*$ is said to be uniformly asymptotically stable in $\widetilde D$.  In particular,  if there exist positive constants $k$ and $ \lambda$, independent of $t_0$, such that 
			\begin{align*}
				\rho(x^*,x(t))\leq k\rho(x^*,x(t_0))e^{-\lambda(t-t_{0})},\quad\forall\; t\geq t_{0}\geq 0 \;\textrm{and} \; x(t_0) \in \widetilde{D},
			\end{align*}
we say $x^*$ is exponentially stable  in the subset $\widetilde{D}$.\end{Definition}	

\begin{Definition}
If $x^*$ is a uniformly asymptotically stable equilibrium point of system (\ref{zuic}), the domain of attraction of $x^*$ is defined as the largest set \(\mathcal{A} \subset M\) such that for every initial time \(t_0 \in [0, +\infty)\) and every initial state \(x_0 \in \mathcal{A}\), the corresponding solution $x(\cdot; t_0, x_0)$ is defined on $[t_0, +\infty)$ and $x(t; t_0, x_0) \to x^*$ as $t \to +\infty$, uniformly in $t_0$.  In particular, if $\mathcal A=M$,  the equilibrium point $x^*$ is said to be globally uniformly 
asymptotically stable.
\end{Definition}

 This paper focuses on the Lyapunov stability theory for System (\ref{zuic}). Specifically, we investigate whether the classical Lyapunov theorem holds for this system and, if so, estimate the {domain of attraction} for its equilibrium point.

In the classical Lyapunov stability theory, the Lyapunov theorem states that if a Lyapunov function can be found for an equilibrium point, then this equilibrium point is {uniformly asymptotically stable}. Moreover, the standard proof of this theorem often provides an explicit {estimate} of the domain of attraction.

The significance of the Lyapunov theorem and the associated estimate of the domain of attraction is twofold:
\begin{enumerate}
    \item {Theoretical Foundation}: The theorem provides a powerful and elegant framework for analyzing stability {without} explicitly solving the system's differential equations. It shifts the problem from solving complex dynamics to finding a suitable scalar function with specific energy-like properties.

    \item {Practical Design \& Safety Guarantee}: The estimated domain of attraction is not merely a theoretical byproduct; it has crucial practical implications. It quantitatively defines a {safe region of operation}. For any initial condition starting within this estimated region, convergence to the equilibrium point  is guaranteed. This is fundamental for controller design, robustness analysis, and certifying system safety in applications such as robotics, power systems, and aerospace engineering. As noted by Khalil in his classic work \cite[Section 8.2, p.312]{jind}, estimating the domain of attraction for a uniformly asymptotically stable equilibrium point is equally important as determining whether an equilibrium point is uniformly asymptotically stable.
\end{enumerate}

	We study the stability of nonautonomous systems on Riemannian manifolds for the following reasons.
\begin{enumerate}
\item Mechanical systems are naturally described by ordinary differential equations on a Riemannian manifold, as discussed in \cite{jix1,jix2,jix3}. The Riemannian metric in this context is determined by the kinetic energy of the system, as shown in \cite[Lemmas 4.29 \& 4.30, p.169]{jix3}. 
When such systems are subject to time-dependent external forces, the resulting equations of motion become nonautonomous ordinary differential equations on the manifold (see \cite[Proposition 4.59, p.195]{jix3}).

\item Existing research on Lyapunov theory for systems on Riemannian manifolds has primarily focused on autonomous systems, including \cite[Section 6.1]{jix3}, \cite[Corollary 20]{MR3139532}, and \cite[Proposition 4.1 \& Corollary 4.5]{dis}, among others.  However, in these results, for uniformly asymptotically stable equilibrium points, either no estimate of the domain of attraction is provided, or it is concluded that the equilibrium point  is globally uniformly asymptotically stable. 
\item  As noted in \cite[Theorem 1]{obs} and \cite[Corollary 5.9.13, p.251]{s}, many systems on manifolds do not possess globally uniformly asymptotically stable equilibrium points. Therefore, from a practical perspective, estimating the domain of attraction of uniformly asymptotically stable equilibrium points becomes crucial for such systems.

\item There is a scarcity of literature on Lyapunov theorems for nonautonomous systems on Riemannian manifolds. Furthermore, existing results do not provide estimates of the domain of attraction for uniformly asymptotically stable equilibrium points, as seen in \cite[Theorem 3.1]{wu}.

\end{enumerate}

   To obtain a Lyapunov-type theorem for system (\ref{zuic}), we employ a metric-dependent mapping that transforms the system (\ref{zuic}) into a system on the tangent space \(T_{x^*}M\) at the equilibrium point \(x^*\). Note that this transformed system is defined on the flat space \(T_{x^*}M\). By applying the Lyapunov theorem for Euclidean spaces (see \cite[Theorem~4.9, p.152]{jind}) to this system, we derive a Lyapunov-type theorem for the original system (\ref{zuic}), as presented in conclusion (i) of Theorem \ref{the1}.

Furthermore, for equilibrium points whose domain of attraction cannot be determined to be the entire manifold \(M\), but which are uniformly asymptotically stable, we adapt the methods used to prove globally uniformly asymptotic stability in Euclidean spaces. This leads to a more refined estimate of the domain of attraction for such equilibrium points, as shown in conclusions (ii) and (iii) of Theorem \ref{the1}.

Subsequently, building on Theorem \ref{the1}, we provide sufficient conditions for exponentially stable equilibrium points and give an estimate of their domain of attraction, as detailed in Theorem \ref{the2}.

Our results exhibit the following novelties:
\begin{enumerate}

\item For systems on Riemannian manifolds, we provide estimates for the domain of attraction of uniformly asymptotically stable equilibrium points. As noted earlier, even for autonomous systems on Riemannian manifolds, the existing literature largely lacks estimates for the domain of attraction when it is not the entire manifold.
\item When the manifold \(M\) is a Euclidean space equipped with the Euclidean metric, our estimates for the domain of attraction of uniformly asymptotically stable equilibrium points (see conclusions (ii) and (iii) of Theorem \ref{the1} and Theorem \ref{the2}) are more refined than the corresponding results in Euclidean spaces (cf. \cite[Theorem~4.9, p.152 \& Theorem 4.10, p.154]{jind}). This is detailed in  item 6 of Remark \ref{rt21}, Remark \ref{the2}, and   item 3 of Remark \ref{r2d4}.
\item The estimates for the domain of attraction obtained in Theorems \ref{the1} $\&$ \ref{the2} depend on the choice of the Riemannian metric. Examples \ref{exa1} $\&$ \ref{exa2} illustrate that for the same system, different Riemannian metrics yield different estimates of the domain of attraction. Therefore, to obtain as refined an estimate as possible, it is necessary to select an appropriate Riemannian metric (see  items  1 $\&$ 2 of Remark \ref{r2d4}).
\item In Theorems \ref{the1} $\&$  \ref{the2}, the estimates obtained for the domain of  attraction of a uniformly asymptotically stable equilibrium point depend on the injectivity radius of the equilibrium point. When the manifold \(M\) is simply connected, the injectivity radius can be bounded using its sectional curvature, as shown in Proposition \ref{cur1}. In Example \ref{exa1}, we apply Proposition \ref{cur1} to estimate the injectivity radius of the equilibrium point.
\end{enumerate}
      
    This paper is organized as follows. In Section \ref{not}, we introduce our main results and present two illustrative examples. Section \ref{proof} is devoted to the proof of our main results. In Section \ref{ap4}, we provide a lemma that will be used in the calculations for Example \ref{exa1}.
	
\setcounter{equation}{0}	
	\section{Main Results}\label{not}
	\def\theequation{2.\arabic{equation}}

Before stating the main results, we introduce some necessary preliminaries.

\subsection{Preliminaries}\label{prel}

  Let \(\nabla\) denote the Levi-Civita connection associated with the Riemannian metric \(g\). For \(x \in M\), denote by \(\exp_x: T_x M \to M\) the exponential map at \(x\) (see \cite[Section 2.1]{cdz}). For \(r > 0\), we also denote by $B(O_x, r)$ the set $ \{ X \in T_x M;\; |X| < r \}$, and by $B_x(r)$ the set $\{y\in M;\,\rho(x,y)<r\}$. The set $B_x(r)$ is called the geodesic ball centered at $x$ with radius $r$ in the metric  $g$.  Let \(i(x)\) represent the injectivity radius at a point \(x \in M\) (see \cite[p.~189]{pth}), which is defined by
\begin{align}\label{dij}
i(x) = \text{the largest }\, r \, \text{for which}\, \exp_x \colon B(0_x, r) \to B_x( r) \text{ is a diffeomorphism} .
\end{align}
It follows from \cite[Lemma 2.1]{cdz} that $i(x)>0$.

We now introduce the curvature tensor of the manifold. 
Let \([X,Y] = XY - YX\) denote the Lie bracket of vector fields \(X\) and \(Y\).  The curvature tensor \(R\) of \((M,g)\) (see \cite[p.~78]{pth}) is the map that assigns to each pair \(X, Y \in TM\) a linear map \(R(X,Y) \colon TM \to TM\) given by
\[
R(X,Y)Z = \nabla_X \nabla_Y Z - \nabla_Y \nabla_X Z - \nabla_{[X,Y]} Z, \qquad \forall Z \in TM.
\]
If \(X(x), Y(x) \in T_{x} M\) are linearly independent vectors at a point \(x \in M\), 
set \(\Pi = \operatorname{span}\{X(x), Y(x)\}\). 
The sectional curvature \(K(x,\Pi)\) of the plane \(\Pi\) is defined by
\[
K(x,\Pi) = \frac{\langle R(X,Y)Y(x), X(x) \rangle}{|X(x)|^2 |Y(x)|^2 - \langle X(x), Y(x) \rangle^2}.
\]

For a differentiable function \(\omega \colon M \to \mathbb{R}\), we denote by \(d\omega\) its differential (see \cite[Definition~2.8, p.~10]{c}). 
For a differentiable function \(W \colon [0, +\infty) \times M \to \mathbb{R}\), 
we write \(W_t(t, x)\) for the partial derivative with respect to \(t\) and \(d_x W(t, x)\) for the partial differential with respect to \(x\).

   \subsection{Statements of Main Results}

    \begin{Definition}
    	$D \subseteq M$ is an open set containing $x^*$. A function $W: D \to \mathbb R$ is said to be positive definite related to $x^*$, if $W(x^*)=0$ and $W(x)>0$, for each $x \in D \backslash \{x^*\}$.
    \end{Definition}
    
    \begin{Definition}
		A continuous function $\alpha:[0,a)\to [0,+\infty)$ with $a>0$ is said to belong to class $\mathcal K$ if it is strictly increasing and $\alpha(0)=0$. Moreover, if $a=+\infty$ and $\lim_{r\to+\infty}\alpha(r)= +\infty$, $\alpha(\cdot)$ is said to belong to $K_\infty$.
    	A continuous function $\beta:[0,a)\times [0,+\infty)\to [0,+\infty)$ is said to belong to class $\mathcal{KL}$, if for each fixed $s$, the mapping $\beta(r,s)$ belongs to class $\mathcal K$ with respect to $r$ and , for each fixed $r$, the mapping $\beta(r,s)$ is decreasing with respect to $s$ and $\beta(r,s)\to 0$ as $s\to +\infty$.
    \end{Definition}
    
    The following assumption is used in the main results.

  \begin{description}
  \item[(A)]The vector field $f:[0,+\infty) \times D \to TM$ is continuous in $t$ and locally Lipschitz in $x$.
  \end{description}
   Our main results are stated as follows:
   
   	\begin{Theorem}\label{the1}
Assume that assumption (A) holds.   	Let $x^{\ast}$ be an equilibrium point of system (\ref{zuic}) and $D \subseteq B_{x^{\ast}}(i(x^{\ast}))$ be an open set containing $x^{\ast}$.  Let $r_{0}=\sup\{r>0;\,   B_{x^{\ast}}(r)\subseteq D\}\in(0, +\infty)\cup\{+\infty\}$.  If there exists a continuously differentiable function $V:[0,+\infty)\times D\rightarrow \mathbb R$ such that
   	\begin{align}\label{kon}			
   		&W_{1}(x)\leq V(t,x)\leq W_{2}(x),\quad \forall\; (t,x)\in [0,+\infty)\times D, \\
   		\label{daos}& V_{t}(t,x)+d_{x}V(t,x)(f(t,x))\leq -W_{3}(x),\quad \forall\; (t,x)\in [0,+\infty)\times D, 
   	\end{align}		
   	where $W_{i}(\cdot)$ 
   	(with $i=1,2,3$) are continuous positive definite functions related to $x^*$ on $D$.
   	Then, the following assertions hold:
 \begin{itemize}
    \item[(i)] For any $r\in(0,r_0)$ and $c\in\left(0,\min\limits_{\rho(x,x^*)=r} W_1(x)\right)$, there  exists a  class $\mathcal{KL}$ function $\beta$ such that
\begin{align}\label{dasl}
\rho(x^*,x(t))\leq\beta(\rho(x^{\ast},x(t_0)),t-t_{0}),\, \forall \,t \geq t_{0}\geq 0,
\end{align}
holds for all $x(t_0)\in\{x\in\text{cl}\,B_{x^*}(r);\,W_2(x)\leq c\}$, where $\text{cl}\,B_{x^*}(r)$ denotes the closure of the set $B_{x^*}(r)$.
\item[(ii)] If 
	\begin{align}\label{rub1}
\lim_{\substack{\rho(x,x^*)\to r_0,
\\ x\in B_{x^*}(r_0)}
}W_1(x)=+\infty,
\end{align}
then, there exists a class $\mathcal{KL}$ function $\beta$ defined on $[0,r_0)\times[0,+\infty)$ such that inequality (\ref{dasl}) holds for all $x(t_0)\in B_{x^*}(r_0)$.
\item[(iii)] If $r_0<+\infty$ and 
\begin{align}\label{ast3}
&\lim_{\substack{\rho(x,x^*)\to r_0
\\ x\in B_{x^*}(r_0)}}\frac{W_1(x)}{r_0^2-\rho^2(x,x^*)}=+\infty,
\\\label{ast4}& d_x\rho^2(x,x^*)\left(f(t,x)\right)\leq 0,\quad \forall\,(t,x)\in[0,+\infty)\times  B_{x^*}(r_0),
\end{align}
then, there exists a class $\mathcal{KL}$ function $\beta$ defined on $[0,r_0)\times[0,+\infty)$ such that inequality (\ref{dasl}) holds for all $x(t_0)\in B_{x^*}(r_0)$.
\end{itemize}

   \end{Theorem}

   	\begin{Remark}\label{rt21}
We refer to Theorem~\ref{the1} as the Lyapunov-type theorem for systems on a Riemannian manifold. Several remarks are in order:

\begin{enumerate}
    \item  For system~\eqref{zuic}, a $C^1$ function $V \colon [0, +\infty) \times D \to \mathbb{R}$ is called a Lyapunov function for the  equilibrium point $x^*$ if there exist continuous positive definite functions  $W_1(\cdot)$, $W_2(\cdot)$, and $W_3(\cdot)$ on the domain $D$ such that  \eqref{kon} and~\eqref{daos} hold. Theorem~\ref{the1} states that if such a Lyapunov function exists, then for any $r\in(0,r_0)$ and $c\in\left(0,\min\limits_{\rho(x,x^*)=r} W_1(x)\right)$, the system is uniformly  asymptotically stable in  $\{x\in\text{cl}\,B_{x^*}(r);\,W_2(x)\leq c\}$.

    \item The conclusions in Theorem~\ref{the1}, namely assertion~(i) and the case $r_0=+\infty$ in assertion~(ii), are extensions of the Euclidean-space result \cite[Theorem~4.9, p.152]{jind} to Riemannian manifolds. 
Unlike the Euclidean case, our results require the domain $D$ to be contained in the normal neighbourhood $B_{x^*}(i(x^*))$. 
The reason for this requirement is that the proof of the Lyapunov theorem in Euclidean space uses the fact that the squared distance function is smooth. 
For $M=\mathbb{R}^n$, this condition holds automatically, whereas for a general Riemannian manifold, applying the same idea to obtain a Lyapunov-type theorem still requires the smoothness of $\rho^2(x^*,\cdot)$, which is guaranteed only inside $B_{x^*}(i(x^*))$.

\item When the manifold $M$ is simply connected, the injectivity radius $i(x^*)$ itself is closely related to the curvature tensor of the manifold. 
Proposition~\ref{cur1} shows how the injectivity radius can be estimated via the sectional curvature of the manifold. 
When $M=\mathbb{R}^n$, the corresponding injectivity radius is infinite, and the results in Theorem~\ref{the1} reduce exactly to the Euclidean case \cite[Theorem~4.9, p.152]{jind}.

\item In assertion~(ii) of Theorem~\ref{the1} with $r_0<+\infty$, if condition (\ref{rub1}) holds, then the domain of attraction of the equilibrium point $x^*$ of system~\eqref{zuic} is $B_{x^*}(r_0)$. 
    This set is larger than $\{x\in\operatorname{cl}B_{x^*}(r);\; W_2(x)\le c\}$ from assertion~(i), where $r\in(0,r_0)$ and $c\in\bigl(0,\min_{\rho(x,x^*)=r}W_1(x)\bigr)$.

\item Assertion (iii) deals with the same case $r_0<+\infty$ but under the additional assumptions (\ref{ast3}) and (\ref{ast4}). 
    Here the domain of attraction is estimated as $B_{x^*}(r_0)$, which is even larger than the one obtained in assertion (i).

\item Even in the Euclidean case $M=\mathbb{R}^n$, assertions (ii) and (iii) of Theorem~\ref{the1} yield a sharper description of the domain of attraction than \cite[Theorem~4.9, p.152 ]{jind}. Thus, our results provide refinements not present in the classical Euclidean theorem.

\end{enumerate}

   \end{Remark}
   
% \begin{Remark}
%As established in \cite[Corollary 20]{MR3139532}, for a dynamical system with multiple equilibrium points on a smooth manifold admitting a Morse-Lyapunov function, the domain of attraction must possess a particular topological property: it must be deformable to the equilibrium set via continuous retraction. This result profoundly reveals the intrinsic connection between dynamical stability and the topological structure of the manifold.
%\end{Remark}
%\begin{Remark}
%For  systems on a connected paracompact manifold, \cite[Theorem 2.2]{MR231409} states that the domain of attraction of an asymptotically stable equilibrium point is diffeomorphic to a Euclidean space.
%
%\cite[Section 6.1]{jix3} discusses the local Lyapunov stability of autonomous systems on differentiable manifolds.
%
% \end{Remark} 
   
   	The sufficient condition of an equilibrium point to be  exponentially stable  in  a set   is stated as follows.
   
   	\begin{Theorem}\label{the2} Assume that assumption (A) holds. 
   	Let $x^{\ast}$ be an equilibrium point of system (\ref{zuic}) and $D \subseteq B_{x^{\ast}}(i(x^{\ast}))$ be an open set containing $x^{\ast}$. Let $r_{0}=\sup\{r>0;\; x \in B_{x^{\ast}}(r)\subseteq D\}$. If there exists a continuously differential function $V:[0,+\infty)\times D\rightarrow \mathbb R$ satisfying
   	\begin{align}\label{zhi1}
   		&k_{1}\rho^{a}(x^{\ast},x)\leq V(t,x)\leq k_{2}\rho^{a}(x^{\ast},x),\quad\forall\; (t,x)\in [0,+\infty)\times D,
   		\\ \label{zhi2} 
   		&V_{t}(t,x)+d_{x}V(t,x)(f(t,x))\leq -k_{3}\rho^{a}(x^{\ast},x),\quad\forall\; (t,x)\in [0,+\infty)\times D ,
   	\end{align}
   	where $k_{i}$ (with $i=1,2,3$) and a are positive constants.  Then, $x^{\ast}$ is exponentially stable   in $ B_{x^{\ast}}((\frac{k_{1}}{k_{2}})^{\frac{1}{a}}r_{0})$, and satisfies
   	\begin{align}\label{es8}
\rho(x(t),x^*)\leq \left( \frac{k_2}{k_1} \right)^{1/a}\rho(x(t_0),x^*)e^{-\frac{k_3}{ak_2}(t-t_0)},\quad\forall\,t\geq t_0\geq 0,\,x(t_0)\in B_{x^{\ast}}\left(\left(\frac{k_{1}}{k_{2}}\right)^{\frac{1}{a}}r_{0}\right).
\end{align}
Here, when $r_0 = +\infty$, $B_{x^{\ast}}\left( \left( \frac{k_1}{k_2} \right)^{\frac{1}{a}} r_0 \right)$ denotes the manifold $M$.
   	\end{Theorem}
   	
   	\begin{Remark}
   	
The result in \cite[Theorem 4]{ni5} studies autonomous systems on Riemannian manifolds and provides sufficient conditions for asymptotic stability of equilibrium points: if the squared distance function serves as a Lyapunov function, then the domain of attraction of the stable equilibrium point includes any region containing no points conjugate to the equilibrium point along geodesics originating from it.

In comparison, Theorems \ref{the1} $\&$ \ref{the2} developed in this work differ in several  aspects:
\begin{enumerate}
    \item It handles the more general case of nonautonomous systems.
    \item In Theorems \ref{the1} $\&$ \ref{the2}, the Lyapunov function is not limited to the squared distance function. Moreover, these two results characterize the domain of attraction in terms of the injectivity radius of the equilibrium point. As noted in item “3” of Remark \ref{rt21}, the injectivity radius can be estimated using the sectional curvature of the Riemannian manifold. Additionally, Theorem \ref{the2} states that when the Lyapunov function additionally satisfies conditions (\ref{zhi1}) and (\ref{zhi2}), the equilibrium point is exponentially asymptotically stable, and the decay rate of solutions is also characterized. This aspect is not addressed in \cite[Theorem 4]{ni5}.
\end{enumerate}
   	\end{Remark}
   	
   	\begin{Remark}\label{r233}
Theorem \ref{the2} can be viewed as an extension of \cite[Theorem 4.10, p.154]{jind} from Euclidean space to Riemannian manifolds. However, even when \( M = \mathbb{R}^n \), Theorem \ref{the2} differs from \cite[Theorem 4.10, p.154]{jind}.  Specifically, in \cite[Theorem 4.10, p.154]{jind}, when \( r_0 < +\infty \), the domain of attraction of the equilibrium point is not characterized; when \( r_0 = +\infty \), the domain of attraction is the entire Euclidean space. In contrast, Theorem \ref{the2} addresses the case where \( r_0 < +\infty \) and provides an estimate for the domain of attraction.
\end{Remark}

In applying Theorems~\ref{the1} and~\ref{the2}, it is often necessary to estimate the injectivity radius of an equilibrium point. The following proposition provides an estimate of the injectivity radius, which is derived from \cite[Corollary 3.3, p. 202]{c}, \cite[Theorem 1.6]{wo}, \cite[Lemma 6.4.7, p. 258]{pth}, and \cite[Proposition 12.9, p. 352]{lee8}.

    \begin{Proposition}\label{cur1}
    	$(M,g)$ is a complete simply connected, $n$-dimensional  Riemannian manifold. The following properties hold:
    	\begin{description}
    		\item [(i)] If there exists positive constants $\delta$, $\sigma$ such that $0<\frac{\delta}{4}<\sigma\leq K(x,\Pi)\leq \delta$ for any $x \in M$ and any two-dimensional linear subspace $\Pi \subset T_{x}M$, then  $\frac{\pi}{\sqrt{\delta}} \leq i(x) \leq \frac{\pi}{\sqrt{\sigma}}$.
    		
    	    \item [(ii)] If $M$ is compact and $K(x,\Pi) \leq \delta \, (\delta >0)$ holds for any $x \in M$ and any two-dimensional linear subspace $\Pi \subset T_{x}M$, then for each $p \in M$,  $i(p) \geq \min\{\frac{\pi}{\sqrt{\delta}}, \frac{1}{2} \cdot (\textrm{length of the shortest geodesic loop based at p}) \}$, where a geodesic loop   is a non-trivial (i.e., non-constant) geodesic that starts and ends at the same point $p$.
    	    
    	     \item [(iii)] If $M$ has nonpositive sectional curvature, then for each $x \in M$, $i(x)=+ \infty$.
    	\end{description}
    	
    \end{Proposition}

  \subsection{Example}\label{exps}
    In this section, we present two examples to illustrate applications of the main results. Before introducing these examples, we need to introduce some notations. Given \(x \in M\),  by the definition of the injectivity radius (\ref{dij}), the map \(\exp_x: B(O_x, i(x)) \to B_x(i(x))\) is a diffeomorphism. Accordingly, we can define the inverse of \(\exp_x\), denoted by \(\exp_x^{-1}: B_x(i(x)) \to B(O_x, i(x))\). Furthermore, for \(y \in B_x(i(x))\), denote the differential of the map \(\exp_x^{-1}\) at \(y\) by \(d(\exp_x^{-1})_y: T_y M \to T_x M\) (see \cite[Definition 2.8, p.10]{c}). By the definition of the differential, this is a linear map.

    \begin{Example}\label{exa1} Set $H^2=\{x=(x_1,x_2)^\top\in\mathbb R^2;\,x_2>0\}$. Assume $d: [0,+\infty)\times H^2\to\mathbb R$ is a function that is continuous in $t\in[0,+\infty)$ and locally Lipschitz continuous in $x\in H^2$. Let $a>0$ and set $A=(0,a)^\top$.  Moreover, there exists a continuous function $k: H^2\setminus\{A\}\to\mathbb R$ satisfying
    \begin{align}\label{ck1}
k(x)>0,\;\forall\,x\in H^2\setminus\{A\}\quad\textrm{and}\quad
\lim_{x\to A}k(x)\sqrt{x_1^2+(x_2-a)^2}=0,
\end{align}
such that $d(t,x)\geq k(x)$ for all $(t,x)\in[0,+\infty)\times (H^2\setminus\{A\})$.
Set $f_1(t,x)=-2d(t,x)x_1 x_2^2$, $f_2(t,x)=d(t,x)x_2(a^2+x_1^2-x_2^2)$, and $f(t,x)=(f_1(t,x), f_2(t,x))^\top$ for $(t,x)\in [0,+\infty)\times H^2$.
    	 Consider the following system 
    	\begin{align}\label{shuang}
    		\left\{\begin{array}{l}
    			\dot x(t)=f(t,x(t)), \quad \text  t \geq t_0\geq 0, 
    			\\
    			x_{2}(t)>0, \quad \forall\; t \geq t_0\geq 0,
    			 \\ 
    			x(t_{0})=x_0,
    		\end{array}\right.
    	\end{align}
    	where $x_{0} \in H^2$.
    	Then, there exists a class $\mathcal{KL}$ function $\beta$ defined on $[0,+\infty)\times [0,+\infty)$ such that
    	 \begin{align}\label{beh1}
\rho\left(x(t; t_0, x_0),A\right)\leq \beta(\rho(x_0, A),t-t_0),\quad \forall\,t\geq t_0\geq 0,\;\forall\,x_0\in H^2,
\end{align}
where $x(\cdot; t_0, x_0))$ is the solution to system (\ref{shuang}) with initial time $t_0$ and initial state $x_0$, and $\rho(\cdot, \cdot)$ is a distance function of $H^2$, given by
\begin{align}\label{dhm}
\rho(x,\hat x)=\ln\frac{\sqrt{(x_{1}-\hat x_1)^{2}+(x_{2}+\hat x_2)^{2}}+\sqrt{(x_{1}-\hat x_1^2)^{2}+(x_{2}-\hat x_2)^{2}}}{\sqrt{(x_{1}-\hat x_1)^{2}+(x_{2}+\hat x_2)^{2}}-\sqrt{(x_{1}-\hat x_1)^{2}+(x_{2}-\hat x_2)^{2}}},
\end{align}
for all $x=(x_1,x_2)^\top, \hat x=(\hat x_1,\hat x_2)^\top\in H^2$.
    	  \end{Example}
   
   {\it Proof.}\; First, the set $H^2$
 is a differentiable manifold endowed with a metric that makes it a Riemannian manifold. We will compute some relevant geometric quantities.
   
    	Indeed, the set  $H^2$ is a differentiable manifold with 
    	 only one coordinate chart $(H^{2},I_{H^2})$,  where $I|_{H^{2}}$ is the restriction of the identity map $I: \mathbb R^{2}\rightarrow \mathbb R^{2}$.  In this coordinate chart, 
    	 let $\{\frac{\partial}{\partial x_1}\Big|_x, \frac{\partial}{\partial x_2}\Big|_x\}$ be the basis for the tangent space at $x\in H^2$.
    	 We endow the manifold $H^{2}$ with  the Riemannian metric $g$:
    	\begin{align}\label{hbm1}
    		g_{ij}(x)=\left\langle \frac{\partial}{\partial x_i}\Big|_x, \frac{\partial}{\partial x_j}\Big|_x\right\rangle= \frac{\delta_{ij}}{x_{2}^{2}},\quad i,j=1,2,\quad x=(x_1,x_2)^\top\in H^2.
    	\end{align}
    	 The Riemannian manifold $(H^{2},g)$ is called the hyperbolic space of dimension $2$. It was proved in \cite[p.162]{c} that $(H^{2},g)$ is a simply connected and complete Riemannian manifold, with constant  sectional curvature $-1$. 
    	
    Let $\rho(x, \hat x)$ denote the distance between $x,\hat x\in H^2$  in the metric $g$. It follows from
    	\cite[Theorem 33.5.1, p. 614]{for}, that $\rho(x,\hat x)$ is expressed by (\ref{dhm}).
According to Proposition \ref{cur1}, we obtain $i(A)=+\infty$. Thus, we can define $\exp_\cdot^{-1}: H^2\to T H^2$. And by
  Lemma \ref{iem1}, we obtain the expression for $\exp_x^{-1}A$:
\begin{align}\label{leie2}
\exp_{x}^{-1}A=-\frac{\rho(x,A)}{\sqrt{x_{1}^{2}+(a+x_{2})^{2}}\sqrt{x_{1}^{2}+(a-x_{2})^{2}}}\left( 2x_1 x_2^2\frac{\partial}{\partial x_1}+x_2(x_2^2-x_1^2-a^2)\frac{\partial}{\partial x_2} \right).
\end{align}

In the coordinate chart $(H^2, I_{H^2})$, the vector field $f$ can be expressed as follows:
\begin{align}\label{lefe}
f(t,x)=\sum_{i=1}^2 f_i(t,x)\frac{\partial}{\partial x_i},\quad (t,x)\in[0,+\infty)\times H^2.
\end{align}

Then, we verify that this system satisfies the conditions of Theorem \ref{the1}, with 
 \( V(t, x) = \rho^2(x, A) \) for \( (t, x) \in [0, +\infty) \times H^2 \), $W_1(\cdot)=W_2(\cdot)=\rho^2(\cdot, A) $, $r_0=+\infty$, 
 and 
 \begin{align*}
W_3(x)=\left\{\begin{array}{ll}\displaystyle\frac{2k(x)\rho(x,A)\left( 4x_1^2 x_2^2+(x_2^2-x_1^2-a^2)^2 \right)}{\sqrt{x_{1}^{2}+(a+x_{2})^{2}}\sqrt{x_{1}^{2}+(a-x_{2})^{2}}},&\forall\,x\in H^2\setminus\{A\};
\\[2mm] 0,& x=A.
\end{array}
\right.
\end{align*}•
Given that \( W_1(x) = W_2(x) = \rho^2(x, A) \),  condition \eqref{kon} holds.
Moreover,  by means of \cite[(5.9)]{dzgJ},  \cite[Lemma 5.2]{dzgJ}, (\ref{hbm1}),  (\ref{lefe}) and (\ref{leie2}), we derive
\begin{align*}
&d_tV(t,x)+d_x V(t,x)(f(t,x))
\\&=\nabla \rho^2(x,A)(f(t,x))
\\&=-2({\exp_x^{-1}A})^\top(f(t,x))
\\&=-2\langle \exp_x^{-1}A,f(t,x)\rangle
\\ &=\frac{-2d(t,x)\rho(x,A)}{\sqrt{x_{1}^{2}+(a+x_{2})^{2}}\sqrt{x_{1}^{2}+(a-x_{2})^{2}}}\left( 4x_1^2 x_2^2+(x_2^2-x_1^2-a^2)^2 \right)
\\&\leq -W_3(x),
\end{align*}•
for all $(t,x)\in [0,+\infty)\times H^2$, where $\nabla \rho^2(x,A)$ denotes the covariant differential (see \cite[p.~124]{kn}) of the function $\rho^2(\cdot,A)$, and $({\exp_x^{-1}A})^\top$ is the dual covector of the tangent vector $\exp_x^{-1}A$ (see (\ref{ddc})).
By setting \( x_1 = r \cos\theta \) and \( x_2 = r \sin\theta + a\) for $(r,\theta)\in[0, +\infty)\times [0,2\pi)$, one can show that
\[
\frac{\rho(x, A)}{\sqrt{x_1^2 + (a + x_2)^2} \sqrt{x_1^2 + (a- x_2)^2}} \left( 4x_1^2 x_2^2 + (x_2^2 - x_1^2 - a^2)^2 \right) = O(x_1^2 + (x_2 - a)^2),  
\]
   as $x \to A$, where $O(c)$ represents a quantity that is of the same order as $c$ in the limit $c\to 0$. It follows from property (\ref{ck1}) that condition
    	 \eqref{daos} holds.

    	Thus, system (\ref{shuang}) satisfies the conditions in Theorem \ref{the1}, and consequently there exists a class $\mathcal{KL}$ function $\beta$ defined on $[0,+\infty)\times [0,+\infty)$ such that (\ref{beh1}) holds.
Thus,  the equilibrium point $A$ is uniformly asymptotically stable in $H^2$. $\Box$

\begin{Example}\label{exa2} Consider the system in Example \ref{exa1}. Assume all the assumptions  there hold. When we endow the manifold $H^2$ with the Euclidean metric. Then, there exists a class $\mathcal {KL}$ function $\hat\beta: [0,a)\times[0,+\infty)\to\mathbb R$ such that
\begin{align}\label{ene1}
|x(t)-A|\leq \hat\beta\left( |x(t_0)-A|, t-t_0 \right),\quad\forall\,t\geq t_0\geq 0,\;,
\end{align}
holds for all $x(t_0)\in \{x\in H^2; \,|x-A|<a\}$, 
where $|\cdot|$ denotes the Euclidean norm in $\mathbb R^2$. The above inequality implies that system (\ref{shuang}) is asymptotically stable in $\{x\in H^2; |x-A|<a\}$.

\end{Example}

{\it Proof.}\; In the Euclidean metric, the distance between $x\in H^2$ and $A$ is $|x-A|$. Thus, 
\begin{align*}
r_0=\sup\left\{r>0;\,\{x\in H^2;\,|x-A|<r\}\subset H^2
\right\}=a.
\end{align*}
Set $ W_1(x)= W_2(x)= V(t,x)=|x-A|^2$ for all $(t,x)\in[0,+\infty)\times H^2$
and
\begin{align*}
W_3(x)=\left\{\begin{array}{ll}
2 k(x) x_2 \left(x_{2}+a\right)\left(x_{1}^{2}+\left(x_{2}-a\right)^{2}\right),& x\in H^2\setminus\{A\};
\\ 0, & x=A.
\end{array}
\right.
\end{align*}
It  follows from (\ref{ck1})  that  function $W_3(\cdot)$ is continuously  positive definite related to $A$.
  Moreover,  relation (\ref{kon}) stands.
By a direct calculation, we obtain
    	\begin{align*}
\sum_{i=1}^2\frac{\partial}{\partial x_i}V(t,x)f_i(t,x)&=-2d(t,x)x_2(x_1^2x_2-a^2 x_2+a^3+x_1^2a+x_2^3-x_2^2 a)
\\&=-2d(t,x) x_{2}\left(x_{2}+a\right)\left(x_{1}^{2}+\left(x_{2}-a\right)^{2}\right)
\\&\leq -W_3(x)
\end{align*}
for all $(t,x)\in [0,+\infty)\times H^2$.
  Moreover,  relations (\ref{ast3}) and (\ref{ast4}) holds. Thus, it follows from (iii) of Theorem \ref{the1} that, there exists a class $\mathcal{KL}$ function $\hat\beta: [0, a)\times[0,+\infty)\to\mathbb R$ such that inequality (\ref{ene1}) holds for all $x(t_0)$ satisfying $|x(t_0)-A|< a$.
$\Box$

    \begin{Remark}\label{r2d4}
   With respect to  Examples \ref{exa1} and \ref{exa2},  we have several remarks.
   \begin{enumerate}
\item   For  system (\ref{shuang}), applying Theorem \ref{the1} with different Riemannian metrics yields different estimates for the domain of attraction.
     If the  metric  (\ref{hbm1}) is used, any geodesic ball $B_A(r)$ centered at $A$ with an arbitrarily large radius $r$ is contained within the space $H^2$, as illustrated in ``Figure (i)''. However, with the Euclidean metric, the maximal radius of a geodesic ball centered at $A$ in $H^2$ is $a$, which results in a more conservative (i.e., smaller) estimate for the domain of attraction; see ``Figure (ii)''.

\item To obtain a more refined estimate of the domain of attraction for uniformly asymptotically stable equilibrium points, we need to employ a suitable Riemannian metric such that the corresponding $r_0$ (defined in Theorem 2.1) is maximized.

\item  As mentioned in Remark \ref{rt21}, assertion~(i) and the case $r_0=+\infty$ in assertion~(ii) of Theorem \ref{the1}, are extensions of the Euclidean-space result \cite[Theorem~4.9, p.152]{jind}.  When applying  \cite[Theorem~4.9, p.152]{jind} to system (\ref{shuang}), we can only obtain that  the equilibrium point $A$ is uniformly asymptotically stable in the set $\{x; |x-A|\leq r\}$, where $r\in(0,a)$. More exactly, there exists a class $\mathcal {KL}$ function $\hat\beta: [0,r]\times [0,+\infty)\to \mathbb R$ such that   inequality (\ref{ene1}) holds for all $x(t_0)\in\{x\in H^2;\,|x-A|\leq r\}$. This set is smaller than the set $\{x\in H^2; \,|x-A|<a\}$ that was obtained in Example \ref{exa2} as an estimate for the domain of attraction.
\end{enumerate}
    \end{Remark}
 
\begin{figure}[H]
    \centering
    \includegraphics[width=0.8\linewidth]{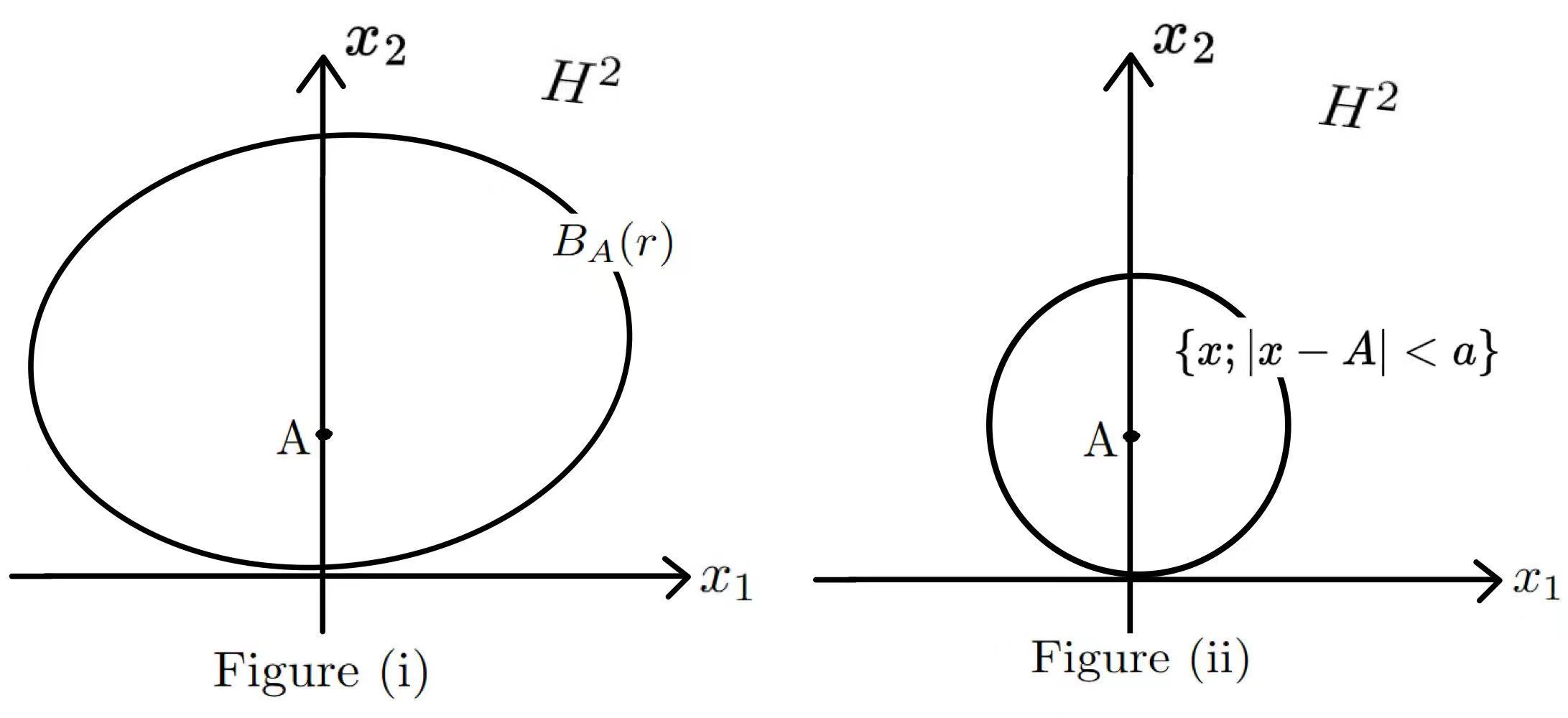}
\end{figure}

	\setcounter{equation}{0}
	\section{Proof of the main results}\label{proof}
\def\theequation{3.\arabic{equation}}
    This section is devoted to proving Theorems~\ref{the1} and~\ref{the2}. 
Before establishing these results, we need two lemmas. 
First, let us introduce some notation. 
For a smooth map $F \colon M \to N$ between two differentiable manifolds $M$ and $N$ and a point $p \in M$, 
we denote by $d(F)_p \colon T_p M \to T_{F(p)} N$ its differential at $p$. 
Our first lemma reads as follows.

	 \begin{lemma}\label{pde}
Let $(M,g)$ be a complete Riemannian manifold and let  $x\in M$. Then, for every $y\in B_x(i(x))$,  there exists a unique ​​minimal geodesic​​ connecting $x$ and $y$, unique up to reparameterization​​. Moreover,  for any $X\in B(O_x,i(x))$, it holds that $\rho(x,\exp_x X)=|X|$.

 \end{lemma}

 {\it Proof.}\;
 First, according to the definition of the injectivity radius (\ref{dij}), we can see that the map $\exp_x: B(O_x, i(x))\to B_x(i(x))$ is a diffeomorphism. For every $y\in B_x(i(x))$, 
 it follows from \cite[Proposition 3.6, p.70]{c} that there exists  a unique minimal geodesic connecting $x$ and $y$,  unique up to reparameterization.
 
 Then, take any $X \in B(O_x,i(x)) $. We obtain from \cite[Proposition 3.6, p.70]{c} that the curve $\gamma: [0,1] \to M$ defined by $\gamma(t) = \exp_x(tX)$ is the minimal geodesic connecting $x$ to $y = \exp_x(X)$. Thus, it follows from the definition of the distance function (see \cite[Definition 2.4, p.146]{c}) and \cite[Corollary 3.9, p.73]{c}  that the distance $\rho(x, y)$ equals the length of $\gamma$, denoted by $\ell(\gamma)$. We refer to \cite[p.43]{c} for the definition of the length of a curve. By the Gauss Lemma (see \cite[Lemma 3.5, p.69]{c}), we have
\begin{align*}
\ell(\gamma) = \int_0^1 |\dot{\gamma}(t)|  dt = \int_0^1 | d(\exp_x)_{tX}(X) |  dt = \int_0^1 |X|  dt = |X|,
\end{align*}
where $d(\exp_x)_{tX}$ denotes the differential of the map $\exp_x$ at the point $tX \in T_x M$.
The proof is complete.
$\Box$
	
\medskip
	
	Our second lemma reads as follows. It is taken from \cite[Proposition~3.6, p.~55]{le1}.

\begin{lemma}\label{l32}
\label{lem:diff_properties_2}
Let $M_1$, $M_2$, and $M_3$ be smooth manifolds with or without boundary, 
let $F_1 \colon M_1 \to M_2$ and $F_2 \colon M_2 \to M_3$ be smooth maps, 
and let $p \in M_1$.
\begin{enumerate}
    \item $d(F_1)_p \colon T_p M_1 \to T_{F_1(p)} M_2$ is linear.
    \item (Chain rule) 
          $d(F_2 \circ F_1)_p = d(F_2)_{F_1(p)} \circ d(F_1)_p 
          \colon T_p M_1 \to T_{F_2 \circ F_1(p)} M_3$.
    \item If $F_1$ is a diffeomorphism, then 
          $d(F_1)_p \colon T_p M_1 \to T_{F_1(p)} M_2$ is an isomorphism, and
  $
          \bigl(d(F_1)_p\bigr)^{-1} = d(F_1^{-1})_{F_1(p)}$.
\end{enumerate}
\end{lemma}

	$\medspace$
	
	\textbf{Proof of Theorem \ref{the1}} We will divide the proof into four steps.
	
{\it Step 1.} We prove assertion (i).

First, we consider a system on the tangent space $T_{x^*}M$. 
Recall the notations for the inverse of $\exp_{x^*}$ and its differential, which are introduced in the first paragraph of Section~\ref{exps}. For each $X \in B(O_{x^*}, i(x^*))$, the differential of the map $\exp_{x^*}^{-1} \colon B_{x^*}(i(x^*)) \to B(O_{x^*}, i(x^*))$ at the point $\exp_{x^*}X \in B_{x^*}(i(x^*))$, denoted by $d(\exp_{x^*}^{-1})_{\exp_{x^*}X} \colon T_{\exp_{x^*}X}M \to T_X(T_{x^*}M)$, is a linear map (see \cite[Definition 2.8, p.10]{c}). Since the tangent space $T_X(T_{x^*}M)$ is isomorphic to $T_{x^*}M$, the map $d(\exp_{x^*}^{-1})_{\exp_{x^*}X}$ can be viewed as a linear map from $T_{\exp_{x^*}X}M$ to $T_{x^*}M$. For any $t \in [0, +\infty)$, since $f(t, \exp_{x^*}X) \in T_{\exp_{x^*}X}M$, it follows that 
$$d(\exp_{x^*}^{-1})_{\exp_{x^*}X} f(t, \exp_{x^*}X) \in T_{x^*}M.$$ Therefore, we can define a map $F \colon [0, +\infty) \times B(O_{x^*}, i(x^*)) \to T_{x^*}M$ by
\begin{align}\label{Ft1}
    F(t, X) = d(\exp_{x^*}^{-1})_{\exp_{x^*} X} f(t, \exp_{x^*}X), \quad \forall\, (t, X) \in [0, +\infty) \times B(O_{x^*}, i(x^*)).
\end{align}
Consequently, we consider the auxiliary system on the tangent space $T_{x^*}M$:
\begin{align}\label{nst}
\dot{X}(t) = F(t, X(t)), \quad \forall\, t \geq t_0 \geq 0.
\end{align}

Second,  for system (\ref{nst}), we construct a Lyapunov function. 
 Define
\begin{align}\label{tvw}\begin{array}{ll}
\tilde{V}(t, X) = V(t, \exp_{x^*}X), & \forall\, (t, X) \in [0, +\infty) \times  \exp_{x^*}^{-1}D, \\
\tilde{W}_i(X) = W_i(\exp_{x^*}X), & \forall\, X \in \exp_{x^*}^{-1}D, \; i = 1, 2, 3.
\end{array}
\end{align}
Since $\exp_{x^*} \colon B(O_{x^*}, i(x^*)) \to B_{x^*}(i(x^*))$ is a diffeomorphism, the function $\tilde{V}$ is continuously differentiable. Note that the tangent space $T_{x^*}M$ is a Euclidean space, endowed with the Riemannian metric $g$ evaluated at $x^*$. We denote the norm of a vector $X \in T_{x^*}M$ by $|X|_{x^*}$.

From \eqref{kon}, we obtain
\begin{align}\label{lyt1}
\tilde{W}_1(X) = W_1(\exp_{x^*}X) \leq V(t, \exp_{x^*}X)=\tilde V(t,X) \leq W_2(\exp_{x^*}X) = \tilde{W}_2(X), 
\end{align}
for all $X \in \exp_{x^*}^{-1}D$.
Let $\tilde{V}_t(t, X)$ and $d_X\tilde{V}(t, X)$ denote the partial differentials of $\tilde V$  with respect to the variables $t$ and   $X$, respectively. For any $X \in B(O_{x^*},i(x^*))$, it follows from Lemma \ref{l32} that $d(\exp_{x^*}^{-1})_{\exp_{x^*}X}$ is the inverse of the linear map $d(\exp_{x^*})_X$. Then, applying the chain rule (see Lemma \ref{l32}) and using \eqref{daos}, we derive
\begin{align*}
&\tilde{V}_t(t, X) + d_X\tilde{V}(t, X)[F(t, X)] \\
&= V_t(t, \exp_{x^*}X) + d_x V(t, \exp_{x^*}X) \circ d(\exp_{x^*})_X \circ d(\exp_{x^*}^{-1})_{\exp_{x^*} X} f(t, \exp_{x^*}X) \\
&= V_t(t, \exp_{x^*}X) + d_x V(t, \exp_{x^*}X) f(t, \exp_{x^*}X) \\
&\leq -W_3(\exp_{x^*}X) = -\tilde{W}_3(X),
\end{align*}
where $(t, X)\in[0,+\infty)\times \exp_{x^*}^{-1}D$.
Thus, $\tilde{V}$ is a Lyapunov function of system \eqref{nst}.

Finally, we  show that, for any $r\in(0,r_0)$ and $c\in\Big(0,\min\limits_{\rho(x,x^*)=r}W_1(x)\Big)$, there exists a  class $\mathcal{KL}$ function $\beta$ such that inequality (\ref{dasl}) holds for all $x(t_0) \in\{x\in\text{cl}\,B_{x^*}(r);\,W_2(x)\leq c\}$.
By means of Lemma 
\ref{pde}, we have
\begin{align*}
\min_{\rho(x,x^*)=r}W_1(x)= \min_{|\exp_{x^*}^{-1}x|_{x^*}=r}\tilde W_1(\exp_{x^*}^{-1}x)=\min_{|X|_{x^*}=r}\tilde W_1(X).
\end{align*}
Thus, $c\in \Big(0,\min_{|X|_{x^*}=r}\tilde W_1(X)\Big)$.

 Applying the Lyapunov theorem for systems on Euclidean space \cite[Theorem 4.9, p.152]{jind}, we conclude that there exists a class $\mathcal{KL}$ function $\beta $ such that
\begin{align}\label{tsbi}
|X(t)|_{x^*} \leq \beta\left(|X_0|_{x^*}, t - t_0\right), \quad \forall\, t \geq t_0 \geq 0, 
\end{align}
for all $X_0 \in\{X\in \text{cl}\, B(O_{x^*}, r);\,\tilde W_2(X)\leq c\}$,
where $X(\cdot)$ is the solution to system \eqref{nst} with initial time $t_0$ and initial state $X_0$. 

Then,  we take any $x_0  \in\{x\in\text{cl}\,B_{x^*}(r);\,W_2(x)\leq c\}$, and let $X(\cdot)$ be the solution to system \eqref{nst} with initial time $t_0 \geq 0$ and initial state $X_0 = \exp_{x^*}^{-1}x_0$. Define $x(t) = \exp_{x^*} X(t)$ for all $t \geq t_0$. Then, by   Lemma \ref{l32} and system \eqref{nst}, we have
\begin{align*}
\dot{x}(t) &= d(\exp_{x^*})_{X(t)} \dot{X}(t) \\
&= d(\exp_{x^*})_{X(t)} F(t, X(t)) \\
&= d(\exp_{x^*})_{X(t)} \circ d(\exp_{x^*}^{-1})_{\exp_{x^*}X(t)} f(t, \exp_{x^*}X(t)) \\
&= f(t, x(t)), \quad \forall\, t \geq t_0.
\end{align*}
Thus, $x(\cdot)$ is the solution to \eqref{zuic} with initial condition $x(t_0) = x_0$. Moreover, from \eqref{tsbi} and Lemma~\ref{pde}, it follows that
\begin{align*}
\rho(x^*, x(t)) = |X(t)|_{x^*} \leq \beta\left(|\exp_{x^*}^{-1}x_0|_{x^*}, t - t_0\right) = \beta\left(\rho(x^*, x_0), t - t_0\right), \quad \forall\, t \geq t_0.
\end{align*}
Since $x_0$ is arbitrary, inequality \eqref{dasl} holds for all initial states starting from $\{x\in\text{cl}\,B_{x^*}(r);\\ W_2(x)\leq c\}$. Thus, we have proved assertion (i) of Theorem \ref{the1}.

{\it Step 2.} We assume  $r_0 =+\infty$,  $D=M$,  and (\ref{rub1}) holds, and show that  there exists a class $\mathcal{KL}$ function such that inequality (\ref{dasl})
 for all $x(t_0)\in M$.

 First, as explained in ``Step 1'', $\tilde V$ is a Lyapunov function for system (\ref{nst}). Moreover, from (\ref{rub1}) and Lemma \ref{pde}, we obtain
\begin{align*}
\lim_{|X|_{x^*}\to +\infty}\tilde W_1(X)
= \lim_{|X|_{x^*}\to +\infty}W_1(\exp_{x^*}X)
= \lim_{\rho(\exp_{x^*}X,x^*)\to +\infty}W_1(\exp_{x^*}X)
= +\infty.
\end{align*}
Applying \cite[Theorem 4.9, p.152]{jind} to system (\ref{nst}), we obtain that there exists a class $\mathcal{KL}$ function $\beta$ such that inequality (\ref{tsbi}) holds for all $X_0 \in T_{x^*}M$. Following the same argument as in ``Step 1'', we conclude that inequality \eqref{dasl} holds for all $x(t_0) \in M$. This proves assertion (ii) for the case $r_0 = +\infty$.

		{\it Step 3.}  
		We prove assertion (ii) for the case that $r_0<+\infty$.

First, we consider a system on $T_{x^*}M$ and show that $O_{x^*}$ is uniformly asymptotically stable in $T_{x^*}M$ for this system. To this end, we define the mapping
		$\varphi: B(O_{x^*}, r_0) \rightarrow T_{x^*}M$ by
\begin{align}\label{vxz}
\begin{array}{l}
\varphi(X) \triangleq (r_0 - |X|_{x^*})^{-1} X,\quad\forall\,X\in  B(O_{x^*}, r_0).
\end{array}
\end{align}
Its inverse mapping $\varphi^{-1}: T_{x^{*}}M \rightarrow B(O_{x^{*}}, r_0)$ is as follows:
\begin{align}\label{ivxz}
\begin{array}{l}
\varphi^{-1}(Y) = (1 + |Y|_{x^*})^{-1} r_0 Y,\quad \forall\,Y\in T_{x^*}M.
\end{array}
\end{align}
We also set the mapping $\hat F: [0,+\infty)\times T_{x^*}M\to T_{x^*}M$ by
\begin{align*}
\hat F(t,Y)=d(\varphi)_{\varphi^{-1}(Y)}F(t,\varphi^{-1}(Y)),\quad \forall\,(t,Y)\in [0,+\infty)\times T_{x^*}M,
\end{align*}
where $F$ is defined by (\ref{Ft1}).
Consider the following system 
\begin{align}\label{svt1}
\dot Y(t)=\hat F(t,Y(t)),\quad\forall\,t\geq t_0\geq 0.
\end{align}
Set
\begin{align*}
\hat V(t,Y)=\tilde V(t,\varphi^{-1}(Y)),\quad \hat W_i(Y)=\tilde W_i(\varphi^{-1}(Y)),\quad\forall\,(t,Y)\in[0,+\infty)\times T_{x^*}M,
\end{align*}
for $i=1,2,3$, where $\tilde V$, $\tilde W_1,\tilde W_2,\tilde W_3$ are defined by (\ref{tvw}). Due to the special structures of the functions $\varphi$ and $\varphi^{-1}$, $\hat W_1(\cdot), \hat W_2(\cdot), \hat W_3(\cdot)$ are continuous positively definite functions related to $O_{x^*}$ on $T_{x^*}M$. Moreover, it follows from Lemma \ref{pde} and (\ref{rub1}) that $\lim_{|Y|_{x^*}\to+\infty}\hat W_1(Y)=+\infty$.

Note that for our case,  (\ref{lyt1}) holds for all $X\in B(O_{x^*}, r_0)$. Thus, we have
\begin{align*}
\hat W_1(Y)=\tilde W_1(\varphi^{-1}(Y))\leq\tilde V(t,\varphi^{-1}(Y))=\hat V(t,Y)\leq \tilde W_2(\varphi^{-1}(Y))=\hat W_2(Y),
\end{align*}
for all $(t,Y)\in[0,+\infty)\times T_{x^*}M$.

Denote by $\hat V_t(t,Y)$ and $d_Y \hat V(t,Y)$ the partial differentials of $\hat V$ with respect to the variables $t$ and $Y$, respectively. Then,  
applying   Lemma \ref{l32}, we derive
\[ 
\begin{array}{l}
d_{Y} \hat{V}(t, Y) [\hat{F}(t, Y)] 
\\
=d_{Y} V\left(t, \exp_{x^*}\varphi^{-1}(Y)\right) \circ d \varphi_{\varphi^{-1}(Y)} F\left(t, \varphi^{-1} (Y)\right)
 \\
=d_{x} V\left(t, \exp _{x^*}{\varphi^{-1}(Y)}\right) \circ d\left(\exp_{x^*}\right)_{\varphi^{-1}(Y)} \circ d (\varphi^{-1})_{Y}\circ d (\varphi)_{\varphi^{-1}(Y)} F\left(t, \varphi^{-1}(Y)\right) 
\\
=d_{x} V\left(t, \exp _{x^*}\varphi^{-1}(Y)\right) \circ d\left(\exp_{x^*}\right)_{\varphi^{-1}(Y)} \circ d\left(\exp_{x^*}^{-1}\right)_{\exp_{x^*}\varphi^{-1}(Y)} f\left(t, \exp _{x^*}\varphi^{-1}(Y)\right) 
\\
=d_{x} V\left(t, \exp _{x^*}\varphi^{-1}(Y)\right) f\left(t, \exp _{x^*}\varphi^{-1}(Y)\right),
\end{array}
\]
for all $(t,Y)\in [0,+\infty)\times T_{x^*}M$.
Thus, we obtain from \eqref{daos} that
\begin{align*}
&\hat V_t(t,Y)+d_Y\hat V(t,Y)[\hat F(t,Y)]
\\&=V_t\left(t,\exp_{x^*}\varphi^{-1}(Y)\right)+d_{x} V\left(t, \exp _{x^*}\varphi^{-1}(Y)\right) f\left(t, \exp _{x^*}\varphi^{-1}(Y)\right)
\\ &\leq -W_3\left(\exp _{x^*}\varphi^{-1}(Y)\right)
\\ &=-\hat W_3(Y),
\end{align*}
for all $Y\in T_{x^*}M$.
Therefore, $\hat V$ is a Lyapunov function of the system (\ref{svt1}). According to the 
Lyapunov theorem for systems on Euclidean space \cite[Theorem 4.9, p.152]{jind}, we conclude that there exists a class $\mathcal{KL}$ function $\hat\beta \colon [0, +\infty) \times [0, +\infty) \to \mathbb{R}$ such that
\begin{align}\label{hbu1}
|Y(t)|_{x^*}\leq \hat\beta\left(|Y(t_0)|_{x^*},t-t_0  \right),\quad \forall\,t\geq t_0\geq 0,\;Y(t_0)\in T_{x^*}M,
\end{align}
where $Y(\cdot)$ is the solution to system (\ref{svt1}).

Then,  
we show that, any solution of (\ref{zuic}) starting from $B_{x^*}(r_0)$ will stay in $B_{x^*}(r_0)$. To this end,  take any $x_0\in B_{x^*}(r_0)$. Then, $\varphi(\exp_{x^*}^{-1}x_0)\in T_{x^*}M$. Denote by $Y_{x_0}(\cdot)$ the solution to system (\ref{svt1}), with initial condition $Y_{x_0}(t_0)=\varphi(\exp_{x^*}^{-1}x_0)$.
 Set 
 \begin{align}\label{vey1}
x(t)=\exp_{x^*}\left( \varphi^{-1}(Y_{x_0}(t)) \right),\quad\forall\,t\geq t_0\geq 0.
\end{align}
 Then, we have 
 \begin{align*}
x(t_0)=\exp_{x^*}\left( \varphi^{-1}(\varphi(\exp_{x^*}^{-1}x_0)) \right)=x_0.
\end{align*}
Moreover, by means of   Lemma \ref{l32} and system (\ref{svt1}), we can derive
\begin{align*}
\dot x(t)&=d\left(\exp_{x^*}  \right)_{\varphi^{-1}(Y_{x_0}(t))}\circ d\left(\varphi^{-1}  \right)_{Y_{x_0}(t)}\dot Y_{x_0}(t)
\\& =d\left(\exp_{x^*}  \right)_{\varphi^{-1}(Y_{x_0}(t))}\circ d\left(\varphi^{-1}  \right)_{Y_{x_0}(t)}\circ  d(\varphi)_{\varphi^{-1}(Y_{x_0}(t))}F\left( t, \varphi^{-1}(Y_{x_0}(t)) \right)
\\&=d\left(\exp_{x^*}  \right)_{\varphi^{-1}(Y_{x_0}(t))}\circ d\left( \exp_{x^*}^{-1} \right)_{\exp_{x^*}\varphi^{-1}(Y_{x_0}(t))}f\left(t, \exp_{x^*}\varphi^{-1}(Y_{x_0}(t))  \right)
\\& =f\left(t, \exp_{x^*}\varphi^{-1}(Y_{x_0}(t))  \right)
\\&=f(t,x(t)),
\end{align*}
for all $t\geq t_0\geq 0$. Thus, $x(\cdot)$ is the solution to the system (\ref{zuic}).
  And  it follows from the definition of $x(\cdot)$ and Lemma \ref{pde} that $x(t)\in B_{x^*}(r_0)$ for all $t\geq t_0$. 

Finally, we show that
there  exists a class $\mathcal{KL}$ function $\beta$ such that inequality  (\ref{dasl}) holds for all $x(t_0)\in B_{x^*}(i(x^*))$.

 To this end, we set
 \begin{align*}
\beta(s,t)=r_0\hat\beta\left(\frac{s}{r_0-s}, t-t_0 \right),\quad \forall\,(s,t)\in[0,r_0)\times [0,+\infty).
\end{align*}
It can be checked that $\beta$ is a class $\mathcal{KL}$ function on $[0,r_0)\times [0,+\infty)$. For any $x_0\in B_{x^*}(r_0)$,
by mean of Lemma \ref{pde}, (\ref{vxz}), (\ref{ivxz}) (\ref{hbu1}) and (\ref{vey1}),  we can derive
\begin{align*}
\rho(x(t), x^*)&=|\exp_{x^*}^{-1}x(t)|_{x^*}
\\&= \left|\exp_{x^*}^{-1}\circ\exp_{x^*}\left(\varphi^{-1}(Y_{x_0}(t))\right)\right|_{x^*}
\\&=\left|\varphi^{-1}(Y_{x_0}(t) )\right|_{x^*}
\\&\leq r_0|Y_{x_0}(t)|_{x^*}
\\&\leq r_0\hat\beta\left(|Y_{x_0}(t_0)|_{x^*}, t-t_0  \right)
\\&= r_0 \hat\beta\left(\left|\frac{\exp_{x^*}^{-1}x_0}{r_0-|\exp_{x^*}^{-1}x_0|_{x^*}}  \right|_{x^*}, t-t_0  \right)
\\&=r_0 \hat\beta\left(\frac{\rho(x^*, x_0)}{r_0-\rho(x_0,x^*)}, t-t_0  \right)
\\&=\beta\left(\rho(x^*,x_0),t-t_0\right),
\end{align*}
where $t\geq t_0\geq 0$.
Thus, $x^*$ is uniformly asymptotically stable in $B_{x^*}(r_0)$. 

{\it Step 4.}\,We now prove assertion (iii). Assume $r_0<+\infty$ and that conditions (\ref{ast3}) and (\ref{ast4}) hold. Define
\[
\bar W_i(x)=\frac{W_i(x)}{r_0^2-\rho^2(x,x^*)}\quad\text{for }x\in B_{x^*}(r_0),\ i=1,2,3.
\]
Then each $\bar W_i(\cdot)$ is a continuous positive definite function with respect to $x^*$ on $B_{x^*}(r_0)$. 
Similarly, set
\[
\bar V(t,x)=\frac{V(t,x)}{r_0^2-\rho^2(x,x^*)}\quad\text{for }(t,x)\in[0,+\infty)\times B_{x^*}(r_0).
\]
From (\ref{kon}) and (\ref{ast3}) we obtain
\begin{align*}
&\bar W_1(x) \le \bar V(t,x) \le \bar W_2(x), 
\qquad \forall\,(t,x)\in[0,+\infty)\times B_{x^*}(r_0),\\[2mm]
&\lim_{\substack{\rho(x,x^*)\to r_0,\\ x\in B_{x^*}(r_0)}} \bar W_1(x) = +\infty.
\end{align*}
Moreover, using (\ref{daos}), (\ref{ast4}) and a direct computation gives
\begin{align*}
&\bar V_t(t,x)+d_x\bar V(t,x)(f(t,x)) \\
&=\bigl(r_0^2-\rho^2(x,x^*)\bigr)^{-1}
\Bigl(V_t(t,x)+d_xV(t,x)(f(t,x))\Bigr) \\
&\qquad +V(t,x)(r_0^2-\rho^2(x,x^*))^{-2}\,d_x\rho^2(x,x^*)(f(t,x)) \\
&\le -\bar W_3(x)
\end{align*}
for all $(t,x)\in[0,+\infty)\times B_{x^*}(r_0)$.

Applying assertion (ii), there exists a class $\mathcal{KL}$ function $\beta$ defined on $[0,r_0)\times[0,+\infty)$ such that inequality (\ref{dasl}) holds for every $x(t_0)\in B_{x^*}(i(x^*))$.

The proof is completed. 
	$\Box$

	$\medspace$

	\textbf{Proof of Theorem \ref{the2}}\; 
	The proof is composed of two steps.

\noindent
\textit{Step 1.}\; We show that, for any $t_0\in[0,+\infty)$ and $x(t_0)\in B_{x^*}\left( \left(\frac{k_1}{k_2}\right)^{1/a} r_0 \right)$, system (\ref{zuic}) admits a global solution $x(\cdot)$. We will prove it in two cases.

\noindent
\textit{Case i}: $r_0=+\infty$. We obtain from (ii) of Theorem \ref{the1} that there exists a class $\mathcal{K}\mathcal{L}$ function $\beta: [0,+\infty)\times[0,+\infty)\to\mathbb{R}$ such that inequality (\ref{dasl}) holds for all $x(t_0)\in B_{x^*}(r_0)=B_{x^*}\left( \left(\frac{k_1}{k_2}\right)^{1/a} r_0 \right)=M$.

\noindent
\textit{Case ii}: $r_0<+\infty$. Take any $t_0\in[0,+\infty)$. Fix any $x(t_0)\in B_{x^*}\left( \left(\frac{k_1}{k_2}\right)^{1/a} r_0 \right)$ and set
\[
r_{x(t_0)}=\frac{1}{2}\left(\frac{k_1}{k_2}\right)^{-1/a}\left(\rho\left(x(t_0), x^*  \right)+\left(\frac{k_1}{k_2}\right)^{1/a}r_0\right)<r_0.
\]
Then,
\begin{align}\label{krt0}
\left(\frac{k_1}{k_2}\right)^{-1/a}\rho(x(t_0),x^*)<r_{x(t_0)}<\left(\frac{k_1}{k_2}\right)^{-1/a} \left(\frac{k_1}{k_2}\right)^{1/a} r_0=r_0.
\end{align}
Set $c=k_1\left(\frac{k_1}{k_2}\right)^{-1}\rho(x(t_0), x^*)^a$. It follows from (\ref{krt0}) that
\[
c<k_1r_{x(t_0)}^a=\min_{\rho(x,x^*)=r_{x(t_0)}}k_1\rho(x,x^*)^a.
\]
By means of (i) of Theorem \ref{the1}, there exists a class $\mathcal{K}\mathcal{L}$ function $\beta_{x(t_0)}: [0,r_{x(t_0)}]\times[0,+\infty)\to\mathbb{R}$ such that
\begin{align}\label{bxtz}
\rho\left(x(t;t_0^\prime,x_0), x^*\right)\leq \beta_{x(t_0)}\left(\rho(x_0, x^*), t-t_0^\prime  \right),\;\forall\,t\geq t_0^\prime\geq 0,
\end{align}
for all $x_0\in \operatorname{cl} B_{x^*}\left(\left( \frac{k_1}{k_2} \right)^{1/a} r_{x(t_0)} \right)$, where $x(\cdot; t_0^\prime,x_0)$ is the solution to system (\ref{zuic}) with initial time $t_0^\prime$ and initial state $x_0$. It follows from (\ref{krt0}) that $x(t_0)\in B_{x^*} \left(\left( \frac{k_1}{k_2} \right)^{1/a} r_{x(t_0)} \right)$. Thus, we obtain from (\ref{bxtz}) that
\begin{align*}
\rho(x(t), x^*)\leq\beta_{x(t_0)}\left( \rho(x(t_0), x^*),t-t_0 \right),\quad\forall\,t\geq t_0.
\end{align*}
Since $x(t_0) \in B_{x^*}\left( \left(\frac{k_1}{k_2}\right)^{1/a} r_0 \right)$ is arbitrary; so for any initial state $x(t_0) \in B_{x^*}\left( \left(\frac{k_1}{k_2}\right)^{1/a} r_0 \right)$, system (\ref{zuic}) has a global solution $x(\cdot)$.

{\it Step 2.}	We prove inequality (\ref{es8}). Since the vector field $f$ is defined on $[0,+\infty)\times D$, we obtain from ``Step 1'' that, for any initial state  $x(t_0)\in B_{x^*}\left( \left(\frac{k_1}{k_2}\right)^{1/a} r_0 \right)$, the corresponding solution $x(\cdot)$ of system  (\ref{zuic}) satisfies $x(t)\in D$ for $t\in[0,+\infty)$.   
		Thus, by means of  (\ref{zhi1}) and (\ref{zhi2}), we obtain
		\begin{align*}
			\frac{d}{dt}V(t,x(t))\leq-\frac{k_{3}}{k_{2}}V(t,x(t)),\quad \forall \, t\geq t_{0}\geq 0 \;\textrm{and}\;x(t_0) \in B_{x^{\ast}}\left(\left(\frac{k_{1}}{k_{2}}\right)^{\frac{1}{a}}r_{0}\right).
		\end{align*}
By means of  Gronwall's inequality, we derive
		\begin{align*}
			V(t,x(t))\leq V(t_{0},x(t_0))e^{-\frac{k_{3}}{k_{2}}(t-t_{0})},\quad \forall \, t\geq t_{0}\geq0\;\textrm{and}\;x(t_0) \in B_{x^{\ast}}\left(\left(\frac{k_{1}}{k_{2}}\right)^{\frac{1}{a}}r_{0}\right).
		\end{align*}
		Recalling (\ref{zhi1}) and (\ref{zhi2}), we obtain from the above inequality that		
\begin{align*}
			\rho(x^{\ast},x(t)) &\leq \Big[\frac{V(t,x(t))}{k_{1}}\Big]^{\frac{1}{a}}
			\\&\leq \Big[\frac{V(t_{0},x(t_0))e^{-\frac{k_{3}}{k_{2}}(t-t_{0})}}{k_{1}}\Big]^{\frac{1}{a}} \nonumber \\
			&\leq\Big[\frac{k_{2}\rho^a(x^{\ast},x(t_0)) e^{-\frac{k_{3}}{k_{2}}(t-t_{0})}}{k_{1}}\Big]^{\frac{1}{a}} \nonumber \\
			&=(\frac{k_{2}}{k_{1}})^{\frac{1}{a}}\rho(x^{\ast},x(t_0))e^{-\frac{k_{3}}{ak_{2}}(t-t_{0})},\quad \forall \; t \geq t_{0}\geq 0,\;\forall \; x(t_0) \in B_{x^{\ast}}\left(\left(\frac{k_{1}}{k_{2}}\right)^{\frac{1}{a}}r_{0}\right).
		\end{align*}			   
Thus, we obtain inequality (\ref{es8}).
	$\Box$

  %%%%%%%%%%%%%%%  

 \setcounter{equation}{0}
 \section{Appdenix}\label{ap4}
 \def\theequation{4.\arabic{equation}}
 
 First, we introduce the notion of the dual (covector) of a tangent vector at a point  $x \in M$, as well as the definition of the dual (vector) of a cotangent vector at $x \in M$, all of which will be used later in Section \ref{proof}.  For a tangent vector $X \in T_{x}M$, we denote by $X^T \in T_{x}^{*}M$ the dual covector of $X$, defined by
 \begin{align}\label{ddc}
{X}^T(Y) = \langle X, Y \rangle, \quad \forall\, Y \in T_{x}M.
\end{align}
Similarly, denote by ${\eta}^T \in T_{x}M$ the dual vector of $\eta \in T_{x}^{*}M$, defined by
\begin{align}\label{ddv}
\langle {\eta}^T, Y \rangle = \eta(Y), \quad \forall\, Y \in T_{x}M.
\end{align}

  The following lemma shows how to compute the inverse map of the exponential map within a coordinate chart.

 	\begin{lemma}\label{iem1}
 		Suppose $M$ is a complete, $n$-dimensional Riemannian manifold. 
Let $x^* \in M$ and let $(U, \varphi)$ be a coordinate chart containing $x^*$, 
with $\varphi(x) = (\xi_1, \dots, \xi_n) \in \mathbb{R}^n$ for all $x \in U$. 
For any $x \in U$, write $\xi = \varphi(x)$ and denote by 
$
\left\{ \left. \frac{\partial}{\partial \xi_1} \right|_x, \dots, \left. \frac{\partial}{\partial \xi_n} \right|_x \right\}
$
the coordinate basis of $T_x M$. 
Let $G(\xi) = \bigl(g_{ij}(\xi)\bigr)_{i,j=1}^n$ be the matrix of the metric tensor with respect to this basis, i.e.,
$
g_{ij}(\xi) = \left\langle \left. \frac{\partial}{\partial \xi_i} \right|_x, \left. \frac{\partial}{\partial \xi_j} \right|_x \right\rangle$ for $ i,j = 1,\dots,n$, 
and let $G^{-1}(\xi) = \bigl(g^{ij}(\xi)\bigr)_{i,j=1}^n$ denote its inverse matrix.
Define the set $A$ as $A = \{ x \in U;\, \rho(x, x^*) < \min\{i(x), i(x^*)\} \}$. Then, for each $x \in A$, we can define $\exp_x^{-1} x^* \in T_x M$. Moreover, the local expression for $\exp_x^{-1} x^*$ is given by
 		 $\exp_{x}^{-1}x^*=-\frac{1}{2}\sum\limits_{i,j=1}^ng^{ij}(\xi)a_{j}(\xi)\left.\frac\partial{\partial\xi_j}\right|_x,$ where $a_{j}(\xi)=\left.\frac\partial{\partial\xi_j}\right|_x\Big(\rho^{2}(\varphi^{-1}(\xi),x^*)\Big)$ and $\varphi^{-1}$ is the inverse map of $\varphi$.
 	\end{lemma}

 	\textbf{Proof: }\; First, according to the definition of injectivity radius (\ref{dij}), we can define $\exp_{x}^{-1}x^*$ for each $x\in A$. Moreover, it follows from \cite[Lemma 5.2]{dzgJ} that
 	\begin{align}\label{bia}
 		\exp_{x}^{-1}x^*=-\frac{1}{2}({\nabla_{1}\rho^{2}(x,x^*)})^\top,\quad \forall \; x \in A,
 	\end{align}
where $\nabla_1\rho^2(x,x^*)$ is the covariant  differential of the function $\rho^2(\cdot,x^*)$, and $({\nabla_{1}\rho^{2}(x,x^*)})^T$ is the dual vector of the cotangent  vector $\nabla_{1}\rho^{2}(x,x^*)$.    We refer to \cite[Definition 5.7, p.102]{c} for the definition of the covariant differential of a tensor field.

Second, we write $\nabla_1\rho^2(x,x^*)$ in the coordinate chart $(U,\varphi)$.
 For any $x\in U$, set $\xi=\varphi(x)$. Then we have
 	$
 		\rho^2(x,x^*)=\rho^2(\varphi^{-1}(\xi),x^*).
$
 	Let $\{d\xi_{1}|_{x}, \dots, d\xi_{n}|_{x}\}$ be the dual basis of $T_{x}^*M$ corresponding to the coordinate basis $\left\{\left.\frac{\partial}{\partial\xi_1}\right|_x, \dots, \left.\frac{\partial}{\partial\xi_n}\right|_x\right\}$, so that
$
d\xi_i|_x\left(\left. \frac{\partial}{\partial\xi_j}\right|_x \right) = \delta_i^j $ for  $i,j=1,\dots,n$,
where $\delta_i^j$ denotes the Kronecker delta. Then, we get 
 $
 		\nabla_{1}\rho^{2}(x,x^*)=\sum\limits_{i=1}^n a_{i}(\xi)d\xi_{i}  |_{x}$, 
 	where 
 	$$a_{i}(\xi)=\nabla_1\rho^2(x,x^*)\left( \left.\frac\partial{\partial\xi_i}\right|_x\right)=\left.\frac\partial{\partial\xi_i}\right|_x(\rho^{2}(\varphi^{-1}(\xi),x^*)),\;i=1,\cdots,n.$$
 	Here the second equality follows from the definition of the covariant differential (see \cite[Definition 5.7, p.102]{c}).
 	
Finally, we express $({\nabla_{1}\rho^{2}(x,x^*)})^T$ in the  coordinate chart.  Assume 
\begin{align*}
({\nabla_{1}\rho^{2}(x,x^*)})^T=\sum\limits_{i=1}^n b_{i}(\xi)\left.\frac\partial{\partial\xi_i}\right|_x.
\end{align*}•
 Recalling the definitions of dual vector and the covariant differential (see \cite[Definition 5.7, p.102]{c}), we have 
 	\begin{align}\label{ai}
 		\left\langle ({\nabla_{1}\rho^{2}(x,x^*)})^T, \left. \frac\partial{\partial\xi_i}\right|_x\right \rangle ={\nabla_{1}\rho^{2}(x,x^*)}\left(   \left. \frac\partial{\partial\xi_i}\right|_x
 		\right)
 		=\left.\frac\partial{\partial\xi_i}\right|_x\Big(\rho^{2}(\varphi^{-1}(\xi),x^*)\Big)=a_{i}(\xi),
 	\end{align}
 	for $i=1,\cdots,n$.
 	On the other hand, applying the expression of $({\nabla_{1}\rho^{2}(x,x^*)})^T$, we get
 	\begin{align}\label{bi}
 		\left\langle ({\nabla_{1}\rho^{2}(x,x^*)})^T, \left. \frac\partial{\partial\xi_i}\right|_x\right \rangle =\left\langle \sum\limits_{j=1}^n b_{j}(\xi)\left.\frac\partial{\partial\xi_j}\right|_x, \left. \frac\partial{\partial\xi_i}\right|_x\right \rangle=\sum\limits_{j=1}^n b_{j}(\xi)g_{ji}(\xi), \quad i=1,\cdots,n.
 	\end{align}
 	Combining (\ref{ai}), (\ref{bi}) and the  symmetry of the matrix $G(\xi)$ (see \cite[p.13]{lee8})​,   we deduce
 	\begin{align}\label{da}
 		\Big(b_{1}(\xi),\cdots,b_{n}(\xi)\Big)^\top=G^{-1}(\xi)\Big(a_{1}(\xi),\cdots,a_{n}(\xi)\Big)^\top.
 	\end{align}

 	According to (\ref{bia}) and (\ref{da}), we obtain the expression of $\exp_{x}^{-1}x^*$ in the local coordinate:
 	\begin{align*}
 		\exp_{x}^{-1}x^*=-\frac{1}{2}\sum\limits_{i,j=1}^ng^{ij}(\xi)\left.\frac\partial{\partial\xi_j}\right|_x\Big(\rho^{2}(\varphi^{-1}(\xi),x^*)\Big) \left.\frac\partial{\partial\xi_i}\right|_x,\quad \forall \; x\in U.
 	\end{align*}
 	The proof is complete.
$\Box$

  \bibliographystyle{amsplain}
\bibliography{D:/studies/tex/deng(3)}
\providecommand{\bysame}{\leavevmode\hbox to3em{\hrulefill}\thinspace}
\providecommand{\MR}{\relax\ifhmode\unskip\space\fi MR }
\providecommand{\MRhref}[2]{%
  \href{http://www.ams.org/mathscinet-getitem?mr=#1}{#2}
}
\providecommand{\href}[2]{#2}

\end{document}